\documentclass{article}
\pdfoutput=1
\usepackage{mathpazo}
\usepackage[utf8]{inputenc}
\usepackage[T1]{fontenc}
\usepackage[english]{babel}
\usepackage[a4paper]{geometry}
\usepackage{amsmath, pifont}
\usepackage[round]{natbib}
\usepackage{hyperref}

\newtheorem{Allitas}{Proposition}[section]
\newtheorem{Def}[Allitas]{Definition}
\newenvironment{Bizonyitas}{\par\noindent{\sc Proof:}\newline}{\hfill\ding{111}}
\DeclareMathOperator{\intt}{int}
\DeclareMathOperator{\cl}{cl}

\author{Gyula Magyarkuti 
\footnote{e-mail: \url{magyarkuti@member.ams.org}
web: \url{http://www.uni-corvinus.hu/magyarkuti }}\\
Corvinus University of Budapest, 
Department of Mathematics\footnote{Address: Budapest, IX Fővám tér 1-3,  (1828 Budapest. P.O. Box 489)}}
\title{Discussion on Lechicki and Spakowski's counterexample}
\date{\today}

\begin{document}
\maketitle
\begin{abstract}
    It is well-known that intersection of continuous correspondences can lost the continuity property.
    Lechicki and Spakowski's theorem says that intersection of H-lsc functions remains H-lsc if the intersection is a bounded subset of a normed space and its interior is nonempty.
    Lechicki and Spakowski pointed to the importance of the boundedness assumption in the case of
    infinite dimensional range giving a counterexample.
    Even though the counterexample works properly and 
    is one of the most cited patterns of discontinuity,
    it has no detailed discussion in the literature of economics and optimization theory.
    What is more, some misleading interpretation of this very important counterexample can be observed.

    Our technical note clarifies the exact role of Lechicki and Spakowski's counterexample,
    computing each of the important properties of the correspondences rigorously. 
\end{abstract}

\section{Introduction}

There are many optimization problems of economic theory, 
when the feasible set is the intersection of moving subsets. 
In virtue of above, the feasible set can be written in the form 
$F\left( x \right)=F_1(x)\cap F_2\left( x \right)$, 
where $x$ is a parameter of the optimization problem.
Discussing these types of problems the very first step is to know the continuity property
of the intersection correspondence.
How does the intersection function $F$ inherit the continuity property of functions $F_1$ and $F_2$?
In economic background 
this problem goes back to \cite{Hildenbrand-1974}.
It is widely discussed in \cite{MR1485775}, using different types of hypertopologies.
Results extending the earlier concepts 
and references to its applications 
can be found in papers 
\cite{springerlink:10.1007/BF00940557}, \cite{springerlink:10.1007/BF00940558}.
\cite{MR1906849} later introduced a unified approach of at least three different types
of continuity of the intersection correspondence and proved an intersection theorem using his new concept.

 \cite{MR796459} gave one of the best and most known solution of this problem.
They show that if $F_1,F_2$ are Hausdorff-lower semi continuous (H-lsc for further reference) convex, closed valued correspondences at $y_0$,
for which $F(y_0)$ is bounded and the interiority assumption 
$\intt F(y_0)\neq\emptyset$ holds, than the intersection correspondence $F$ remains H-lsc at $y_0$.
It was also proved that the boundedness assumption is superfluous in the case when the range set of correspondences is finite dimensional normed space.
They gave a beautiful counterexample where the range is infinite dimensional normed space 
$\ell_\infty$ and the set $F\left( y_0 \right)$ is unbounded. \cite[Example 3.]{MR796459}.

Since 1985 this counterexample has become one of the most cited pattern of discontinuity of the infinite dimensional intersection problem. 
See for example \citet{MR1485775}, \citet{MR1676280}, and \cite{MR1189795}.
We have never seen
a detailed discussion of this important example, not even in handbooks.
The book \citet{MR1485775} devoted to be a foundation of multivalued analysis misinterprets this 
counterexample. 
According to that the intersection correspondence fails to be Vietoris-lower semi continuous (V-lsc). 
Actually, $F$ is really not H-lsc, but it does have the V-lsc property.

The goal of our paper is to clarify the role of Lechicki and Spakowski's counterexample 
computing each of the important properties of the functions rigorously.
Even if some of them can really be considered as a classroom case, 
some of them are far from trivial.

\section{Preliminary notes}
Let $\mathcal{P}(X)$ denotes the power set of a normed space $X$ and 
$B_r$ denotes the open ball centered to zero, with radiuos $r$ as usual. 
The normed space of bounded sequences endowed with the supremum norm is denoted by $\ell_\infty$.
The vector $t\in\ell_\infty$ denotes $t=\left( t_1,t_2,\dots,t_n,\dots \right)$.
The $k$-th unit vector multiplied by $k$ is denoted by 
$e^{\left( k \right)}=\left( 0,\dots,k,\dots \right)$.

Recall the basic definitions of the continuity of correspondences.
    Suppose that $S$ is a topological space.
    Consider a closed-valued correspondence $F$ of $S$ into $X$.
    The correspondence 
    is called to be \emph{lower hemicontinuous} or \emph{Vietoris-lower continuous, (V-lsc)}
    at the element $y\in S$,
    if for every open set $G\subseteq X$ such that $F\left( y \right)\cap G\neq\emptyset$
    there exists a neighborhood $V$ of $y$ for which 
    $F(y')\cap G\neq\emptyset$ for all $y'\in V$.\\
    A harder concept is the \emph{Hausdorff-lower semi continuity}.
    The function $F$ above is \emph{H-lsc} at $y$, 
    for every $r>0$, there exists a neighborhood $V$ of $y$ for which
    $
    F\left( y \right)\subseteq F\left( y' \right)+B_r
    $
    for all $y'\in V$.

The well-known relationship is that H-lsc implies V-lsc and the two concepts are the same
if the values are totally bounded subsets of $X$. 
See \cite[e.g.]{Michael-1951}.

The following real functions play important role later.
If parameters $x\geq 0$ and $t_1\in\mathbb{R}$ are given, consider the functions
$f_1,f_2:\left\{ 2,3,\cdots \right\}\to \mathbb{R}$ 
\[
f_1\left( k \right)= k\left( 1-t_1-x \right)\quad\text{and}\quad
f_2\left( k \right)= k+\frac{1}{k}\left( t_1-x \right).
\]
Easy computation shows that
\[
f_2(k)-f_1(k)=
\frac{1}{k}(t_1-x) + k(x+t_1).
\]
Thus $f_2(k)\leq f_1(k)$ if and only if 
\[
k^2\left( x+t_1 \right)\leq x - t_1.\tag{\dag}
\]

Now, we should separate the cases with respect to $t_1$ and $x$.\\
        First, $t_1\geq 0$.
        At that time $x-t_1\leq x\leq x+t_1$ since ($\dag$) has no solution with $k\geq 2$,
        except $t_1=0=x$.
        It means that 
        \[
        f_1(k)\leq f_2(k), \quad \forall k\geq 2, k\in\mathbb{N}
        \]\\
        Second, $t_1\leq -x$.
        At that time $t_1 + x \leq 0$ thus the left hand side of ($\dag$) is 
        $k^2\left( x+t_1 \right)\leq 0$
         and
        for the right hand side of ($\dag$) is $x-t_1\geq x\geq 0$ because $t_1$ is non positive.
        It means that
        \[
        f_2(k)\leq f_1(k), \quad \forall k\geq 2, k\in\mathbb{N}
        \]\\
        The third case $-x<t_1<0$ is never used.

\section{The Spakowski Lechicki's example}
\subsection{Simple part}
We introduce a correspondence $F_1$ 
which is closed, convex valued and H-lsc at 0.
\begin{Def}
Define the correspondence $F_1:[0,1]\to\mathcal{P}\left( \ell_\infty \right)$ is as follows.
\[
F_1(x)=
\left\{ t\in \ell_\infty:
t_1\geq x, 
t_k\leq k-x \text{ if }
k\geq 2
\right\}
\]
\end{Def}
One can easily verify that $F_1\left( x \right) $ is a convex, closed subset of the normed space $\ell_\infty$.
The important special case at zero is
\[
F_1(0)=
\left\{ t\in \ell_\infty:
t_1\geq 0, 
t_k\leq k \text{ if }
k\geq 2
\right\}
\]
The H-lsc property of $F_1$ is our first step.
\begin{Allitas}
    For every $\varepsilon>0$ the inequality $0\leq x< \varepsilon/2$ implies
    $F_1\left( 0 \right)\subseteq F_1\left( x \right)+B_\varepsilon$.
    It means that $F_1$ is H-lsc at zero.
\end{Allitas}

\begin{Bizonyitas}
If $s\in F_1\left( 0 \right)$, then $s_1\geq 0, s_k\leq k$ that is
$s_1+\varepsilon\geq \varepsilon/2>x$ and for $k\geq 2$ the $s_k-\varepsilon/2\leq k-\varepsilon/2< k-x$.
Introduce a vector belonging to $F_1\left( x \right)$ as
$t=\left( s_1+\varepsilon/2 ,s_2-\varepsilon/2 ,s_3-\varepsilon/2 ,\cdots \right)$.
It is clear that 
$\|t-s\|=\varepsilon/2$, thus $s\in t+B_\varepsilon$. It was to be proved.
\end{Bizonyitas}
\subsection{Tweak part}
The second function is more interesting. 
\begin{Def}
Define the correspondence $F_2:[0,1]\to \mathcal{P}\left( \ell_\infty \right)$ as
\[
F_2(x)=
\left\{ 
t\in \ell_\infty:
t_1 \leq 1-x,
t_k\leq \min\left\{
k\left( 1-t_1-x \right),
k+\frac{t_1}{k}-\frac{x}{k}
\right\}
\forall k\geq 2
\right\}
\]
\end{Def}
As a corollary of the discussion in section 2, we obtain the special case at zero.
\[
F_2(0)=
\left\{ 
t\in \ell_\infty:
t_1 \leq 1,
\text{ but if } k\geq 2\text{ then }
\begin{cases}
    t_k\leq k\left( 1-t_1 \right)&\text{, if } t_1\geq 0\\
    t_k\leq k+\frac{t_1}{k}&\text{, if } t_1<0
\end{cases}
\right\}
\]
\begin{Allitas}
    The set $F_2\left( x \right)\subseteq \ell_\infty$ is closed, convex set for every $x\in\left[ 0,1 \right]$.
\end{Allitas}
\begin{Bizonyitas}
    Consider the sequences $s,t\in F_2\left( x \right)$ and for $\lambda\in\left[ 0,1 \right]$ denote
    $r=\lambda s + \left( 1-\lambda \right)t$.
    For the first term $r_1\leq\max\left\{ s_1,t_1 \right\}\leq 1-x$. 
    Suppose that $k\geq 2$ is an integer and $x$ is given.
    The function
    $
    g\left( t \right)=\min\left\{ k\left( 1-t-x \right),k+\frac{t}{k}-\frac{x}{k} \right\} 
    $
    is a concave function, because the minima of concave functions remains concave.
    Thus $r_k=\lambda s_k+\left( 1-\lambda \right)t_k\leq
    \lambda g\left( s_1 \right)+\left( 1-\lambda \right)g\left( t_1 \right)\leq
    g\left( \lambda s_1+\left( 1-\lambda \right)t_1 \right)=
    g\left( r_1 \right)
    $ which means the convexity of $F_2\left( x \right)$.
    Similarly, $F_2\left( x \right)$ is closed by the continuity of $g$ defined above.
\end{Bizonyitas}
\begin{Allitas}
    Using the notations above,
    for every $\varepsilon>0$ the inequality $0\leq x< \varepsilon /2$ implies
    $
    F_2\left( 0 \right)\subseteq F_2\left( x \right) + B_{\varepsilon}.
    $
    Thus $F_2$ is also H-lsc at zero.
\end{Allitas}

\begin{Bizonyitas}
    Suppose that $0<x<\varepsilon /2$ and $s\in F_2\left( 0 \right)$  are given.

    If $s_1>\varepsilon /2$, then $s_1-x> 0$. Thus introducing $t_1=s_1-x$
    and $t_k=s_k \forall k\geq 2$ we obtain that $t\in\ell_\infty$ and $\|t-s\|=s_1-t_1=x<\varepsilon/2$.
    But $t_k=s_k\leq k\left( 1-s_1 \right)=k\left( 1-\left( s_1-x \right)-x \right)=
    k\left( 1-t_1-x \right)$ which proves that $t\in F_2\left( x \right)$.
    Summarizing this case of the proof the inclusion $s\in F_2\left( x \right)+B_{\varepsilon/2}$ holds true.

    If $s_1<-\varepsilon$, then $s_1<-2x$. Thus if $t_1$ denotes $t_1=s_1+x$ we see that $t_1<-x$.
    Let us introduce $t_k=s_k$ for all $k \geq 2$. 
    Clearly $t\in\ell_\infty$ and $\|s-t\|=t_1-s_1=x<\varepsilon/2$ holds.
    But $t_k=s_k\leq k+\frac{s_1}{k}=k+\frac{s_1+x}{k}-\frac{x}{k}=k+\frac{t_1}{k}-\frac{x}{k}$,
    thus $t\in F_2\left( x \right)$.
    As in the previous case the inclusion $s\in F_2\left( x \right)+B_{\varepsilon/2}$ also holds.

    Finally suppose $-\varepsilon\leq s_1\leq\varepsilon /2$. 
    Let $t_1=-x$. 
    Introduce $t_k=\min\left\{ s_k,k-\frac{2x}{k} \right\}$ for all $k\geq 2$. 
    The fact $s\in\ell_\infty$ assures $t\in\ell_\infty$.
    In virtue of the negativity of $t_1$ the inclusion $t\in F_2\left( x \right)$ holds true.
    Realize that $\frac{2x}{k}\leq x < \varepsilon /2$ for all $k\geq 2$.
    Thus in the case when $t_k\neq s_k$ we obtain
    \[
    t_k< s_k\leq k< k-\frac{2x}{k}+\frac{\varepsilon}{2}=t_k+\frac{\varepsilon}{2}.
    \]
    So $|s_k-t_k|<\frac{\varepsilon}{2}$ for $k\geq 2$ and clearly 
    $|s_1-t_1|<\varepsilon$.
    Thus the inequality $\|t-s\|<\varepsilon$ must hold, 
    so $s\in F_2\left( x \right)+B_{\varepsilon}$ holds true in the third case as well.
\end{Bizonyitas}

\subsection{Intersection}
Define the intersection of $F_1$ and $F_2$.
\begin{Allitas}
    Let $F\left( x \right)=F_1\left( x \right)\cap F_2\left( x \right)$ be the intersection correspondence
    for every $x\in\left[ 0,1 \right]$. Then
    \[
    F\left( x \right)=
    \left\{ 
    t\in \ell_\infty:
    x\leq t_1\leq 1-x,
    t_k\leq k\left( 1-t_1-x \right)
    \forall k\geq 2
    \right\}
    \]
    is a closed, convex valued correspondence.
\end{Allitas}

The special case at zero is:
\[
F\left( 0 \right)=
\left\{ 
t\in \ell_\infty:
0\leq t_1\leq 1,
t_k\leq k\left( 1-t_1 \right) \forall k\geq 2
\right\}
\]

\begin{Allitas}
    Denote the set $S=\left\{ e^{\left( k \right)}:k\in \mathbb{N} \right\}$.
    Then for every $x>0$ and for every $r>0$ the inclusion $S\subseteq F\left( x \right)+B_{r}$ does not hold.
    Considering to $S\subseteq F\left( 0 \right)$ the correspondence $F$ is not H-lsc at zero.
\end{Allitas}
\begin{Bizonyitas}
    Suppose contrary, that there exist $r>0$ and $x>0$, such that 
    for every $k$ there exists $t^{\left( k \right)}\in F\left( x \right)$, 
    with $\|e^{\left( k \right)}-t^{\left( k \right)}\|_{l_\infty}<r$.
    Then the inequality 
    $k<t_k^{\left( k \right)} + r$ holds true for every $k\geq 2$.
    But
    $t_k^{\left( k \right)}\leq k\left( 1-\left( t_1^{\left( k \right)} + x \right) \right)
    \leq k\left( 1-x \right)<k-r$, 
    for for fixed $x$ and for some $k$ which is big enough, but it is a contradiction.
\end{Bizonyitas}

We mention that if 
$t\in F\left( 0 \right), t=\left( t_1,0,\dots \right)$ then $t+B_\delta\subseteq F\left( 0 \right)$ if
$\delta=\min\left\{ \frac{2}{3}\left( 1-t_1 \right),t_1 \right\}$.
Thus the interior of $F\left( 0 \right)$ is nonempty.

The following proposition shows that the weaker continuity concept does not fail at this case.
To prove the V-lsc property,
we could use some intersection theorem as well for example \cite[Lemma 2.6]{MR1301514} or 
\cite{BarabashBusygin-1978}\nocite{BerlyandMarkin-1978} but a direct proof is preferred here.
This is the misinterpreted Example~2.53 of \cite{MR1485775}.
\begin{Allitas}
    The correspondence $F$ is V-lsc at $0$.
\end{Allitas}
\begin{Bizonyitas}
    Let us fix the numbers $t\in F\left( 0 \right)$ and $r>0$. 
    We are going to show, that there exists a number $\delta>0$ such that
    for every $0\leq x <\delta$ the 
    $F\left( x \right)\cap\left( t+\cl\left( B_r \right) \right)\neq\emptyset$ holds true.
    
    The easier case is when $t_1>0$.
    Define the number $\delta$ with $\delta<\min\left\{ t_1/2,r \right\}$.
    If $x<\delta$ we introduce the sequence $s\in\ell_\infty$ as  $s_1=t_1-x$ 
    and $s_k=t_k$ for every $k\geq 2$. 
    Of course $s_1>t_1-t_1/2=t_1/2\geq x$ 
    and $s_1\leq 1-x$.
    It is also clear that $s_k\leq k\left( 1-t_1 \right)=k\left( 1-s_1-x \right)$ that is 
    $s\in F\left( x \right)$.
    In virtue of definition of $\delta$ we obtain that $\|s-t\|_{l_\infty}=t_1-s_1=x<r$.
    Thus our sequence $s$ has the property $s\in F\left( x \right)\cap\left( t+B_r \right)$.

    Consider the case $t_1=0$.
    Choose the number $M$ so big that $M>\|t\|+1$ and $\frac{r}{2M}<1/2$ both hold.
    Denote $\delta=\frac{r}{2M}$.
    If $x<\delta$ is given we define the sequence $s\in\ell_\infty$ as follows 
    \[
    s_k=\left\{
    \begin{array}{ll}
        \delta &,\mbox{if }k=1;\\
        k\left( 1-2\delta \right)&, \mbox{if } 2\leq k\leq M \mbox{ and } k-r\leq t_k;\\
        t_k-r &,\mbox{otherwise.}
    \end{array}
    \right.
    \]
    It is clear that $s\in \ell_\infty$ and $x\leq s_1\leq 1-x$.
    Now, we are going to show that 
    \[
    s_k\leq k\left( 1-2\delta \right), \forall k\geq 2 \tag{1}
    \]
    If $2\leq k\leq M$ then 
    $k-r=
    k\left( 1-\frac{r}{k} \right)\leq
    k\left( 1-\frac{r}{M} \right)=
    k \left( 1-2\delta \right)<k$.
    Thus if $t_k<k-r$ then inequality $(1)$ above really holds.\\
    If $k>M$ then 
    $s_k=t_k-r<M-r=M\left( 1-\frac{r}{M} \right)<k\left( 1-2\delta \right)$
    also holds true. \\
    Thus we can summarize, that $s_k\leq k\left( 1-2\delta \right)\leq k\left( 1-s_1-x \right)$, because
    $s_1+x=\delta + x \leq 2\delta$. It means that inclusion $s\in F\left( x \right)$ holds.
    
    Finally, compute the distance $\|s-t\|$.
    Clearly, $s_1-t_1=\delta\leq \frac{r}{2M}<r$.
    It is easy to see that inequality $|s_k-t_k|\leq r$ holds for all $2\leq k\leq M$.
    Surely, if $k\geq M$ then $t_k-s_k=r$, thus 
    $\|s-t\|=r$.
    We proved, that $F\left( x \right)\cap\left(t+ \cl\left( B_r \right) \right)\neq\emptyset$.
\end{Bizonyitas}

\bibliographystyle{plainnat}
\bibliography{SL}
\end{document}